\documentclass{amsart}

\newcommand{\Z}{\ensuremath{\mathbf Z}}
\newcommand{\Zd}{\ensuremath{ {\mathbf Z}^d }}

\newtheorem{theorem}{Theorem}
\newtheorem{lemma}{Lemma}
\newtheorem{corollary}{Corollary}
\newcommand{\bt}{\begin{theorem}}
\newcommand{\et}{\end{theorem}}
\newcommand{\bl}{\begin{lemma}}
\newcommand{\el}{\end{lemma}}
\newcommand{\bc}{\begin{corollary}}
\newcommand{\ec}{\end{corollary}}
\newcommand{\beq}{\begin{equation}}
\newcommand{\eeq}{\end{equation}}
\newcommand{\benum}{\begin{enumerate}}
\newcommand{\eenum}{\end{enumerate}}

\begin{document}

\title{Sets with more sums than differences}

\author{Melvyn B. Nathanson}
\address{Lehman College (CUNY),Bronx, New York 10468}
\email{melvyn.nathanson@lehman.cuny.edu}
\date{\today}

\thanks{This is the text of a lecture at the SIAM Conference on Discrete Mathematics in Victoria, British Columbia, on June 27, 2006.  This work was supported in part by grants from the NSA Mathematical Sciences Program and the PSC-CUNY Research Award Program.}

\maketitle

\begin{abstract}
Let $A$ be a finite subset of the integers or, more generally, of any abelian group, written additively.  The set $A$ has \emph{more sums than differences} if $|A+A|>|A-A|.$  A set with this property is called an \emph{MSTD set}.
This paper gives explicit constructions of families of MSTD sets of integers.
\end{abstract}

\section{MSTD sets}
Let $G$ be an abelian group, written additively, and let $A$ and $B$ be  subsets of $G$.    We define the {\em sumset}
\[
A+B = \{a+b : a \in A, b\in B\}
\]
and the {\em difference set}
\[
A-B = \{a-b : a \in A, b\in B\}.
\]
In particular, if $A=B$, then 
\[
A+A = \{a+a' : a,a' \in A\}
\]
and 
\[
A-A = \{a-a' : a,a' \in A\}.
\]

Let $[a,b]$ denote the interval of integers $\{a,a+1,a+2,\ldots,b\}$, and let $|A|$ denote the cardinality of the set $A$.

Since $G$ is abelian, for every pair $a, b$ of distinct elements of $G$ we have
\[
a+b=b+a
\]
but
\[
a-b\neq b-a
\]
unless $a-b$ has order 2.
This suggests that if $A$ is finite, then we should ``usually'' have $|A-A| \geq |A+A|$.  Thus, sets with {\it more sums than differences  (MSTD)} should be rare.  On the other hand, elements are omitted from $A-A$ in pairs: If $c \notin A-A$, then $-c \notin A-A$.  This provides an opportunity to construct sets with more sums than differences.
We shall call a finite set $A$ with $|A+A|>|A-A|$ an {\it MSTD set}.

For example, the set
\beq  \label{example1}
A_1 = \{0,2,3,4,7,11,12,14\}
\eeq
is an MSTD set of integers with $A_1+A_1 = [0,28]\setminus\{1,20,27\}$ and $A_1-A_1 = [-14,14]\setminus \{\pm 6,\pm 13\}$, hence $|A_1+A_1| = 26$ and $|A_1-A_1| = 25$.
The set
\beq  \label{example2}
A_2 = \{0,2,3,4,7,9,13,14,16\}
\eeq
is also an MSTD set, with $A_2+A_2 = [0,32]\setminus\{1,24,31\}$ and $A_2-A_2 = [-16,16]\setminus \{\pm 8,\pm 15\}$, hence $|A_2+A_2| = 30$ and $|A_2-A_2| = 29.$

For every set $A$ of integers and every integer $x$, we define the {\it dilation} $x*A = \{xa : a \in A\}.$  If $A$ is an MSTD set of integers, then, for all integers $x\neq 0$ and $y$, the set
\[
x*A+\{y\} = \{xa+y : a\in A\}
\]
is also an MSTD set, and so the property of having MSTD is an affine invariant of the ring of integers.

In this paper we provide explicit constructions of infinite families of MSTD sets that contains the sets $A_1$ and $A_2$, and we also give an algorithm to construct MSTD sets of integers from MSTD sets in finite abelian groups.
This answers a question posed by Nathanson~\cite{nath06x}.

The first examples of MSTD sets of integers appear in Marica~\cite{mari69} and Freiman-Pigarev~\cite{frei-piga73}.   Ruzsa~\cite{ruzs79,ruzs84,ruzs92c} used probability methods to show the existence of MSTD sets.
Roesler~\cite{roes00} proved that, in the collection of all sets of size $k$ in $\{1,2,\ldots,n\}$, the average number of sums does not exceed the average number of differences.  On the other hand, O'Bryant~\cite{obry06} recently proved that a positive proportion of the sets in  $\{1,2,\ldots,n\}$ have more sums than differences.

\section{A construction of MSTD sets of integers}

The subset $A$ of an abelian group $G$ is called {\em symmetric with respect to $a^*$} if $a^*\in G$ and
\[
A = a^*-A.
\]
If $A$ is symmetric with respect to $a^*$, then $A$ has the same number of sums as differences, since
\[
|A+A| = |A+(a^*-A)| = |a^*+(A-A)| = |A-A|.
\]
For example, the arithmetic progression $\{a,a+m,a+2m,\ldots,a+(k-1)m\}$ 
is symmetric with respect to $a^* = 2a+(k-1)m$.  For every $d \geq 1,$ the $d$-dimensional generalized arithmetic progression
\[
L = \{a+x_1m_1+\cdots + x_dm_d : \ell_i \leq x_i \leq \ell_i + k_i-1 \text{ for } i = 1,\ldots,d\}
\]
in the group $G$ is symmetric with respect to
\[
a^*  = 2a+   \sum_{i=1}^d \left( 2\ell_i + k_i-1)m_i\right).
\]
If $A$ is a finite set of integers that is symmetric with respect to $a^*$, then $a^* = \min(A)+\max(A).$

In this section we shall construct infinite families of MSTD sets of integers.  Each MSTD set will be constructed by adjoining a single element to a symmetric set that is a small perturbation of a generalized arithmetic progression.

We begin with an examination of the set 
$A_1  = \{0,2,3,4,7,11,12,14\}$.   
Observe that  $A_1\setminus \{4\}$ is symmetric with respect to 14, and that
\begin{align*}
A_1 & = \{0,2,3,4,7,11,12,14\} \\
& = \{0,2,3,7,11,12,14\} \cup \{ 4\} \\
& = \{0,2\} \cup \{3,7,11\} \cup \left( 14 - \{0,2\}\right)  \cup \{ 4\}.
\end{align*}
Hegarty and Roesler (personal communication) observed that 
this set belongs to the following infinite sequence of MSTD sets:
For $k \geq 3$, let
\[
A^* = \{ 0,2\} \cup \{3,7,11,\ldots, 4k-1\} \cup \{ 4k,4k+2\}
\]
and 
\[
A = A^* \cup \{ 4\}.
\]
To see that the set $A$ is an MSTD set of integers, we note that$A^*$ is symmetric and so $|A^*+A^*| = |A^*-A^*|$.  Since $8 \in (A+A) \setminus (A^*+A^*)$, it suffices to prove that $|A-A| = |A^*-A^*|$.  For this, we only need to show that $A^*-\{ 4\} \subseteq A^*-A^*$, and this is true because 
\begin{align*}
1 & = 3-2 \in A^*-A^* \\
4 & = 7-3 \in A^*-A^* \\
4k-4 & = (4k-1)-3 \in A^*-A^* \\
4k-2 & = (4k)-2 \in A^*-A^*.
\end{align*}

We can generalize this construction as follows.

\bt         \label{MSTD:theorem:construction1}
Let $m, d,$ and $k$ be integers such that $m \geq 4$, $1 \leq d \leq m-1$, $d \neq m/2$, and $k \geq 3$ if $d<m/2$ and $k \geq 4$ if $d > m/2$.
Let
\begin{align*}
B & = [0,m-1] \setminus \{d\} \\
L & = \{m-d,2m-d,\ldots, km-d\} \\
a^* &= (k+1)m-2d \\
A^* &= B \cup L \cup \left( a^* - B\right)\\
A & = A^* \cup \{m\}.
\end{align*}
The set $A$ is an MSTD set of integers.
\et

\begin{proof}
The set $A^*$ is symmetric with respect to $a^*$, and so 
\[
|A^*+A^*| = |A^*-A^*|.
\]
Also, $2m \in A+A$.  We shall prove that $2m \notin A^*+A^*.$

Assume that $2m = a+a' \in A^*+A^*$ with $a \leq a'$.  Then $m \leq a' \leq 2m$ and so $a' \notin B$.  Also, 
\[
\min(a^*-B) \geq a^* - (m-1) = km-2d+1 > 2m
\]
if $(k-2)m \geq  2d$, and this inequality holds  if $k \geq 4$ or if $k \geq 3$ and $d \leq m/2$.  Thus, $a' \notin a^*-B$.  

Since $3m-d>2m$, we have $a' \in L$ only if $a' = 2m-d.$  Therefore, if $a' \in A^*$, then $a' = 2m-d$ and $a=d$.  Since $d \notin B$ and $m-d \neq d$, it follows that $a\notin A^*$ and so $2m \notin A^*+A^*.$  This proves that
\[
2m \in \left(A+A\right) \setminus \left(A^*+A^*\right)
\]
and so
\[
|A+A| > |A^*+A^*|.
\]

We shall prove that $A-A = A^*-A^*.$  
It suffices to show that $A^* - \{ m\} \subseteq A^*-A^*.$  
Note that $0 \in A^*$ implies that $A^* \subseteq A^* - A^*.$

Since $m = (2m-d) - (m-d) \in L-L$, it follows that
\[
\{ m\}  - B = \{1,2,\ldots,m-1,m\}\setminus \{ m-d\} \subseteq B \cup (L-L) \subseteq A^* - A^*
\]
and also $B - \{ m\} \subseteq A^* - A^*.$
Similarly, 
\[
L- \{ m\}  =  \{-d,m-d,\ldots, (k-1)m-d\} \subseteq L \cup \{-d\}.
\]
If $1 < d < m-1$, then $\{1,d+1\}\subseteq B\cup L \subseteq A^*$ and $d \in A^*-A^*.$
If $d=1$, then $\{2,3\} \subseteq B\cup L$ and $1\in A^*-A^*.$
If $d=m-1,$ then $m+1 = 2m-d\in L$, $2 \in B$, and so $m-1 \in B-L   \subseteq A^*-A^*$.  Therefore,
\[
L-\{ m\}  \subseteq A^*-A^*.
\]

Finally, we must prove that
\[
a^*-B- \{ m\}  = \{a^*-b-m : b \in B\} \subseteq A^*-A^*.
\]
If $d = m-1$, then $a^* = (k-1)m+2$ and 
\[
A^* = [0,m-2] \cup \{1,m+1,2m+1,\ldots, (k-1)m+1\} \cup (a^*-[0,m-2]).
\] 
If $b \in [0,m-4]$, then
\[
a^*-b-m = a^*-(b+2)-(m-2) \in A^*-A^*.
\]
If $b = m-3$, then
\[
a^*-b-m = a^*-(m-4)-(m+1) \in A^*-A^*.
\]
If $b = m-2$, then
\[
a^*-b-m = a^*-(m-3)-(m+1) \in A^*-A^*.
\]

Suppose that $1 \leq d \leq m-2.$   
For  $b\notin \{d-1, m-1\}$, we have $b+1 \in B $ and so $a^*-(b+1) \in a^*-B $.  Since $m-1 \in B$, we have
\[
a^*-b-m = \left( a^* - (b+1) \right) - (m-1) \in  A^*-A^*.
\]
If $b = m-1$, then 
\begin{align*}
a^*-b-m & = (k-1)m-2d+1 \\
& = ((k-1)m-d)-(d-1) \\
& \in L-B \subseteq A^*-A^*.
\end{align*}
If $b = d-1$ and $d \neq m-2$, then $2 \leq d \leq m-3$ and
\begin{align*}
a^*-b-m & = a^* - (d-1)-m \\
& = \left( a^* - (d+1) \right) - (m-2) \\
& \in A^* - A^*.
\end{align*}
The last case is $b = d-1$ and $d = m-2$.  Then $m\geq 5$, $k \geq 4$,  and
\[
a^* = (k+1)m-2d  = (k-1)m+4
\]
and
\begin{align*}
a^*-b-m & = (k-1)m+4 - (d-1) - m\\
& = (k-3)m+7\\
& = (k-2)m+2-(m-5) \\
& \in L-B  \subseteq A^*-A^*.
\end{align*}
Therefore, 
\[
a^*-B-m \subseteq A^*-A^*
\]
and so $A-A=A^*-A^*.$
This completes the proof.
\end{proof}

The MSTD set $A_1$ is the case $m=4$, $d=1$, and $k=3$
of Theorem~\ref{MSTD:theorem:construction1}.

Next we consider the MSTD set $A_2 =  \{0,2,3,4,7,9,13,14,16\}$.  This set is a perturbation of a two-dimensional generalized arithmetic progression:
\begin{align*}
A_2 & =  \{0,2,3,4,7,9,13,14,16\} \\
& = \left( \{0,2\} \cup \{3,7\} \cup \{9,13\} \cup \{14,16\} \right) \cup \{4\} \\
& = \{0,2\} \cup \{3+6x_1+4x_2 : 0 \leq x_1 \leq 1 \text{ and } 0 \leq x_2 \leq 1\} \cup \left(16 - \{0,2\}\right) \cup \{4\}.
\end{align*}
The following sequence of sets extends this example.

\bt   \label{MSTD:theorem:example2}
For $k \geq 2$, let
\begin{align*}
A^* & = \{0,2\} \cup \{3,7,11,\ldots, 4k-1\} \cup \{9,13,17,\ldots, 4k+5\} \cup \{4k+6,4k+8\}  \\
& = \{0,2\} \cup \{3+6x_1+4x_2 : 0 \leq x_1 \leq 1 \text{ and } 0 \leq x_2 \leq k-1\} \cup \left(16 - \{0,2\}\right) 
\end{align*}
and
\[
A = A^* \cup \{ 4\}.
\] 
Then $A$ is an MSTD set of integers.
\et

\begin{proof}
It suffices to show that $4-A^* \subseteq A^* -A^*$, and this follows from
\begin{align*}
1 & = 3-2 \in A^*-A^* \\
4 & = 7-3 \in A^*-A^* \\
5 & = 7-2 \in A^*-A^* \\
4k+2 & = (4k+5)-3 \in A^*-A^* \\
4k+4 & = (4k+6)-2 \in A^*-A^*.
\end{align*}
This completes the proof.
\end{proof}

\bl   \label{MSTD:lemma}
Let $m \geq 4$ and let $r$ and $s$ be integers such that 
\[
r \geq 1 \text{ and }  r+1 \leq s \leq m-1.
\] 
Let
\[
B = [0,r-1] \cup [s,m-1].
\]
If 
\beq       \label{MSTD:lemma:2B}
s \leq 2r-1 \text{ and } 2s \leq m+r-1
\eeq
then 
\[
2B = [0,2m-2].
\]
If 
\beq       \label{MSTD:lemma:B-B}
s \leq 2r-1 \text{ or } 2s \leq m+r-1
\eeq
then 
\[
B-B = [1-m,m-1].
\]
Let $m \geq 4$ and $1 \leq r \leq m-2$ and let
\[
B = [0,m-1] \setminus \{r \}.
\]
Then
\[
B-B  = [-(m-1),m-1].
\]
If $2 \leq r \leq m-3$, then  
\[
B+B = [0,2m-2].
\]
If $r=1$, then 
\[
B+B = [0,2m-2]\setminus \{1\}.
\]
If $r=m-2$, then 
\[
B+B =  [0,2m-2] \setminus \{2m-3\}.
\]
\el

\begin{proof}
This is a straightforward calculation.  If $s \leq 2r-1$ and $2s \leq m+r-1,$ then
\[
2A = [0,2r-2] \cup [s,m+d-2] \cup [2s,2m-2] = [0,2m-2].
\]
If $s \leq 2r-1$ or $2s \leq m+r-1,$ then
\[
(A-A) \cap [0,m-1] = [0,r-1] \cup [s-r+1,m-1] \cup [0,m-1-s] = [0,m-1]
\]
and so $A-A = [1-m,m-1]$.  

If $s = r+1$, then $1 \leq r \leq m-1$ and $B = [0,m-1]\setminus \{r\}$.   Condition~\eqref{MSTD:lemma:B-B} always satisfied, and condition~\eqref{MSTD:lemma:2B} is equivalent to $2 \leq r \leq m-3$.  The proofs for the cases $r=1$ and $r=m-2$ are similar.
\end{proof}

The following theorem produces explicit families of MSDT sets constructed from 2-dimensional arithmetic progressions.

\bt   \label{MSTD:theorem:construction2}
Let $m$ and $d$ be integers such that 
\benum
\item[(i)]
$ m\geq 4$
\item[(ii)]
$1 \leq d \leq m-2$ and $d \neq m/2$
\item[(iii)]
if $d < m/2$, then $d \neq m/3$
\item[(iv)]
if $d > m/2$, then $d \neq 2m/3$.  
\eenum
For $k \geq 3$, let
\[
a^* = (k+2)m
\]
\[
B = [0, m-1] \setminus \{ d\}
\]
\[
L = \{2m-d,3m-d,\ldots, km-d\} \cup \{2m+d,3m+d,\ldots,km+d \}
\]
\[
A^* = B \cup L \cup \left( \{a^*\} - B \right)
\]
\[
A = A^* \cup \{m\}.
\]
Then $A$ is an MSTD set of integers.
\et

\begin{proof}
The set $A^*$ is symmetric with respect to $a^*$ and so
\[
|A^*+A^*| = |A^*-A^*|.
\]
Since 
\[
A^* \cap [0,2m] = B \cup \{2m-d\}
\]
and $d \notin B$, we have $2m \notin A^* + A^*$ and so  $2m \in (A+A)\setminus (A^*+A^*)$.  Thus, it suffices to show that $A-A = A^* - A^*$.  We need only prove that
\[
A^* - \{m\} \subseteq A^*-A^*.
\]
Since $d \in B-B$ by Lemma~\ref{MSTD:lemma} and since $m=(3m-d)-(2m-d) \in L-L$, it follows that
\[
\{m\}-B = [1,m]\setminus \{m-d\} \subseteq (B-B) \cup (L-L) \subseteq A^*-A^*
\]
and so
\[
B - \{m\} \subseteq A^*-A^*.
\]

The inequalities $1 \leq d < m/2$ and $d \neq m/3$ imply that $m-2d \neq d$, $m-2d \in B$,  and so
\[
m+d = (2m-d) - (m-2d) \in L-B.
\]
The inequalities $m/2 < d \leq m-2$ and $d \neq 2m/3$ imply that $2m-2d \neq d$, $2m-2d \in B$,   and so
\[
m+d = (3m-d) - (2m-2d) \in L-B. 
\]
Also, $m-d \in B.$  Therefore,
\[
L-\{m\} \subseteq B \cup  L \cup ( L-B) \subseteq A^* - A^*.
\]

Suppose that $1 \leq d \leq m-3$.  If $b \in B$ and $b \neq m-1$, then Lemma~\ref{MSTD:lemma} implies there exist $b_1,b_2 \in B $ such that $m+b = b_1 + b_2$, and so
\[
a^* -b - m = (k+1)m - b = ((k+2)m-b_1) - b_2 \in \left( \{a^*\} - B \right) - B .
\]
If $d = m-2$ and $b = m-3$, then $2m-3 \notin B+B$ but $m\geq 5$ and 
\[
a^*-b-m = km+3 = (km+d)-(d-3) \in (a^* - B) - B.
\]
For all $d \in [1,m-2]$ we have
\[
a^* - (m-1)-m = km+1 = (km+d)-(d-1) \in L - B
\]
and so
\[
 \{a^*\} - B  - \{m\} \subseteq ( \left( \{a^*\} - B \right) - B)  \cup (L - B) \subseteq A^*-A^*.
\]
This completes the proof.
\end {proof}

The next theorem gives a method to construct MSTD sets of integers from generalized arithmetic progressions of dimension $d$ for every $d \geq 1$.

\bt  \label{MSTD:theorem:gap}
Let $m \geq 4$ and let $B$ be a subset of $[0,m-1]$ such that 
\[
B+B = [0,2m-2]
\]
and
\[
B-B = [ -m+1,m-1].
\]
Let $L^*$ be a $(d-1)$-dimensional arithmetic progression contained in $[0,m-1]\setminus B$ such that $\min(L^*)-1\in B$.  Let $k \geq 2$ and let $L$ be the $d$-dimensional arithmetic progression
\[
L = (m-L^*) + m\ast[1,k].
\] 
Let
\[
a^* = \min(L)+\max(L) = (k+3)m-\min(L^*)-\max(L^*)
\]
and
\[
A^* = B \cup L \cup (\{a^*\} - B).
\]
Then $A = A^* \cup \{m\}$ is an MSTD set of integers.
\et

\begin{proof}
The set $A^*$ is symmetric with respect to $a^*$ and so
\[
|A^*+A^*| = |A^*-A^*|.
\]
Since $B \cap L^* = \emptyset$ and
\[
A^* \cap [0,2m] = B  \cup \left( \{2m\} - L^* \right)
\]
we have $2m \notin A^* + A^*$ and so  $2m \in (A+A)\setminus (A^*+A^*)$.  The set $A$ will be an MSTD set if $A-A = A^* - A^*$, and for this it suffices to show that
\[
A^* - \{m\} \subseteq A^*-A^*.
\]

Note that the condition $B+B = [0,2m-2]$ implies that $0,m-1 \in B$.  In particular, $A^* \subseteq A^* - A^*.$

Since $m \in L-L$, we have
\[
\{m\}-B \subseteq [1,m] \subseteq (B-B) \cup (L-L) \subseteq A^*-A^*
\]
and so
\[
B - \{m\} \subseteq A^*-A^*.
\]
Also,
\begin{align*}
L- \{m\} & = (\{m\}-L^*) \cup \left( (m-L^*) + m\ast [1,k-1] \right) \\
& \subseteq (B-B) \cup L \\
& \subseteq A^*-A^*.
\end{align*}
Finally, let $a^*-b -m \in a^* - B - \{m\} .$  
If $b\neq m-1$, then there exist $b_1,b_2 \in B$ such that $m+b = b_1+b_2$, and so
\[
a^*-b -m = (a^* - b_1) - b_2 \in (a^*-B) - B \in A^*-A^*.
\]
If $b = m-1$, then  
\begin{align*}
a^*-b-m & = a^* - 2m+1 \\
& = (k+3)m-\min(L^*)-\max(L^*) -2m+1\\
& = (k+1)m - \max(L^*) - (\min(L^*) - 1) \\
& \in L-B \\
& \subseteq A^* - A^*
\end{align*}
because $(m-\max(L^*)) + km \in L$ and $\min(L^*)-1 \in B$.
This completes the proof.
\end{proof}

Theorem~\ref{MSTD:theorem:gap} is a powerful tool for constructing MSTD sets.  For example, let $P$ be any $(d-1)$-dimensional arithmetic progression of nonnegative integers with $\min(P) = 0$ and $\max(P) = M.$
Choose integers $r$ and $s$ so that
\[
r \geq M+2
\]
and 
\[
r+M+1 \leq s \leq 2r-1.
\]
Choose an integer $m$ such that
\[
2s \leq m+r-1.
\]
Let
\[
B = [0,r-1] \cup [s,m-1]
\]
and
\[
L^* = \{r\}+P.
\]
Then $L^*$ is a $(d-1)$-dimensional arithmetic progression contained in $[0,m-1] \setminus B$, and $\min(L^*)-1 = r-1 \in B$.  
By Lemma~\ref{MSTD:lemma}, we have $B+B = [0,m-2]$ and $B-B = [1-m,m-1].$  
The sets $B$ and $L^*$ satisfy the conditions of Theorem~\ref{MSTD:theorem:gap}.

A slight variation of the proof of Theorem~\ref{MSTD:theorem:gap} gives the following result.

\bt  \label{MSTD:theorem:gap-2}
Let $m \geq 4$ and let $B$ be a subset of $[0,m-1]$ such that 
\[
B+B = [0,2m-2]
\]
and
\[
B-B = [ -m+1,m-1].
\]
Let $L^*$ be a $(d-1)$-dimensional arithmetic progression contained in $[0,m-1]\setminus B$ such that $\min(L^*)-1\in B$ 
and $m \notin L^*+L^*$.  Let $k \geq 2$ and let $L$ be the $d$-dimensional arithmetic progression
\[
L = (m-L^*) + m\ast[0,k].
\] 
Let
\[
a^* = \min(L)+\max(L) = (k+3)m-\min(L^*)-\max(L^*)
\]
and
\[
A^* = B \cup L \cup (\{a^*\} - B).
\]
Then $A = A^* \cup \{m\}$ is an MSTD set of integers.
\et

\section{Finite abelian groups and lattices}
Let $A$ be a nonempty subset of an abelian group.  For all integers $h \geq 0$, we define the {\it $h$-fold sumset}  $hA$ inductively:
\begin{align*}
0A & = \{0\} \\
1A & = A \\
2A & = A+A \\
hA & = (h-1)A+A \text{ for all $h \geq 3$.}
\end{align*}
For nonnegative integers $h$ and $k$ we define the {\it generalized sum-difference set}
\[
\begin{split}
hA-kA = \{ & a_1+\cdots + a_h - a'_1-\cdots - a'_k : a_j \in A \\
& \text{ for $j=1.\ldots,h$ and } a'_{\ell} \in A \text{ for $\ell = 1,\ldots, k$} \}.
\end{split}
\]

Let $G$ be a finite abelian group.  Since every finite abelian group is isomorphic to a direct product of cyclic groups, we can assume that 
\[
G = \Z/m_1\Z \times \cdots \times \Z/m_d\Z,
\]
where $m_1,\ldots,m_d$ are integers greater than 1.  Then
\[
|G| = m_1\cdots m_d.
\]
Let \Zd\ denote the standard $d$-dimensional integer lattice.  The sublattice $\Lambda_G$ of \Zd\ determined by the group $G$ is
\[
\Lambda_G = \{(q_1m_1,\ldots, q_dm_d) \in \Zd : q_i \in \Z \text{ for } i=1,\ldots,d\}.
\]
The fundamental integer parallelepiped in \Zd\ defined by $G$ is the set
\[
P_G = \{ (r_1,\ldots,r_d) \in \Zd : 0 \leq r_i \leq m_i -1 \text{ for } i = 1.\ldots,d\}
\]
Every lattice point in \Zd\ is uniquely the sum of an element in $P_G$ and an element in $\Lambda_G$, that is, 
\beq   \label{MSTD:directsum}
\Zd =  P_G \oplus \Lambda_G .
\eeq
We define the {\em canonical embedding $\varphi$} of $G$ into the fundamental integer parallelepiped $P_G$ of the lattice $\Lambda_G$
by
\beq  \label{MSTD:map1}
\varphi(u_1+m_1\Z, \ldots,u_d+m_d\Z) = (r_1,\ldots,r_d)
\eeq
where $r_i$ is the unique integer such that
\beq  \label{MSTD:map2}
r_i \equiv u_i \pmod{m_i}
\eeq
and
\beq  \label{MSTD:map3}
0 \leq r_i \leq m_i -1.
\eeq
We define a map
\[
\pi: \Zd \rightarrow P_G
\]
by 
\[
(x_1,\ldots,x_d) \mapsto (r_1,\ldots,r_d)
\]
where
\[
r_i \equiv x_i \pmod{m_i}
\]
and
\[ 
0 \leq r_i \leq m_i -1.
\]

For integers $s \leq t$, we define
\[
\Lambda_G(s,t) = \{(q_1m_1,\ldots, q_dm_d) \in \Lambda_G :  s \leq q_i < t \text{ for } i=1,\ldots,d\}.
\]
Then
\beq       \label{MSTD:sizeLambda}
|\Lambda_G(s,t)| = (t-s)^d
\eeq
\beq       \label{MSTD:sumLambda}
\Lambda_G(s_1,t_1) + \Lambda_G(s_2,t_2) = \Lambda_G(s_1+s_2,t_1+t_2-1)
\eeq
and
\beq       \label{MSTD:diffLambda}
\Lambda_G(s_1,t_1) - \Lambda_G(s_2,t_2) = \Lambda_G(s_1-t_2+1,t_1-s_2)
\eeq
In particular, we have the generalized sum-difference relation
\beq       \label{MSTD:sumdiffLambda}
h\Lambda_G(0,t) - k\Lambda_G(0,t) = \Lambda_G(-kt+k,ht-h+1).
\eeq

\bl  \label{MSTD:lemma:lat}
 Let $A$ be a nonempty subset of the finite abelian group $G = \Z/m_1\Z \times \cdots \times \Z/m_d\Z$.  For all integers $h \geq 1$ and $k \geq 0$,
\benum
\item[(i)]
\[
\varphi(hA-kA) = \pi(h\varphi(A)-k\varphi(A)) 
\]
\item[(ii)]
\[
h\varphi(A) - k\varphi(A) \subseteq 
\varphi(hA-kA)  +  \Lambda_G(-k,h).
\]
\item[(iii)]
\[
\varphi(hA-kA) \subseteq 
h\varphi(A) -  k\varphi(A) + \Lambda_G(-h+1,k+1). 
\]
\eenum
\el

\begin{proof}
(i)  If $(s_1,\ldots,s_d) \in \varphi(hA-kA),$ then the group $G$ contains elements $a_1,\ldots,a_h, a'_1,\ldots, a'_k$ such that
\[
(s_1,\ldots,s_d) = \varphi(a_1+\cdots a_h-a'_1-\cdots - a'_k).
\]
For $j=1,\ldots,h$ and $ i=1,\ldots,d$, we can write
\beq \label{MSTD:lat1}
a_j = (r_{1j} + m_1\Z,\ldots, r_{dj} + m_d\Z) 
\eeq
with
\beq \label{MSTD:lat2}
0 \leq r_{ij} \leq m_i-1.
\eeq
Similarly, for $\ell=1,\ldots,k$ and $ i=1,\ldots,d$, we have
\beq \label{MSTD:lat3}
a'_{\ell} = (r'_{1\ell} + m_1\Z,\ldots, r'_{d\ell} + m_d\Z) 
\eeq
with
\beq \label{MSTD:lat4}
0 \leq r'_{i\ell} \leq m_i-1.
\eeq
There exist unique integers $q_i$ and $s'_i$ such that $0 \leq s'_i \leq m_i-1$
and 
\[
\sum_{j=1}^hr_{ij} - \sum_{\ell =1}^k r'_{i\ell} =q_im_i + s'_i.
\]
The $i$th component of the group element
\[
a_1+\cdots + a_h -a'_1-\cdots - a'_{\ell} \in hA-kA
\]
is
\[
\left( \sum_{j=1}^hr_{ij} - \sum_{\ell =1}^k r'_{i\ell}\right) + m_i\Z 
= s'_i + m_i\Z
\]
and so
\[
(s_1,\ldots,s_d) = \varphi\left(  a_1+\cdots + a_h -a'_1-\cdots - a'_{\ell} \right) = (s'_1,\ldots,s'_d).
\]
Therefore, $s_i=s'_i$ for $i=1,\ldots,d$ and 
\beq \label{MSTD:lat5}
\sum_{j=1}^hr_{ij} - \sum_{\ell =1}^k r'_{i\ell} =q_im_i + s_i.
\eeq

Since
\[
\varphi(a_j) = (r_{1j},\ldots,r_{dj}) \in \varphi(A) 
\text{ for $j=1,\ldots, h$}
\]
and
\[
\varphi(a'_{\ell}) = (r'_{1\ell},\ldots,r'_{d\ell})  \in \varphi(A) \text{ for $\ell=1,\ldots, k$}
\]
it follows that
\[ 
\begin{split}
\varphi(a_1) + \cdots & + \varphi(a_h) -  \varphi(a'_1) - \cdots - \varphi(a'_k) \\
& = \left(  \sum_{j=1}^hr_{1j} - \sum_{\ell =1}^k r'_{1\ell}, \ldots, \sum_{j=1}^hr_{dj} - \sum_{\ell =1}^k r'_{d\ell} \right) \\
& = (q_1m_1+s_1,\ldots,q_dm_d+s_d)  \\
& \in h\varphi(A) - k\varphi(A)
 \end{split}
\]
and so
\begin{align*}
\pi\left( \varphi(a_1) + \cdots  + \varphi(a_h) -  \varphi(a'_1) - \cdots - \varphi(a'_k)  \right) 
& = \pi(q_1m_1+s_1,\ldots,q_dm_d+s_d) \\
&= (s_1,\ldots,s_d) \\
& \in \pi(h\varphi(A)-k\varphi(A)).
\end{align*}
Therefore,
\[
\varphi(hA-kA) \subseteq \pi(h\varphi(A)-k\varphi(A)).
\]
Similarly,
\[
\pi(h\varphi(A)-k\varphi(A)) \subseteq \varphi(hA-kA).
\]
This proves~(i).

(ii)  Every lattice point in the generalized sum-difference set $h\varphi(A)-k\varphi(A)$ is of the form 
\[
\varphi(a_1) + \cdots + \varphi(a_h) - \varphi(a'_1) - \cdots - \varphi(a'_k)
\]
with $a_1,\ldots,a_h,a'_1,\ldots,a'_k \in G$.
Using the same notation as above, we see that the $i$th coordinate of this lattice point is 
\beq   \label{MSTD:lat6}
\sum_{j=1}^hr_{ij} - \sum_{\ell =1}^k r'_{i\ell} = q_im_i + s_i
\eeq
where
\[
0 \leq s_i \leq m_i-1
\]
and
\[
(s_1,\ldots,s_d) \in \varphi(hA-kA).
\]
Moreover,
\[
-km_i \leq -k(m_i-1) \leq q_im_i+s_i \leq h(m_i-1) < hm_i
\]
and so
\beq \label{MSTD:lat7}
-k \leq q_i < h
\eeq
for $i=1,\ldots, d$.
Therefore, 
\[
(q_1m_1,\ldots,q_dm_d) \in \Lambda_G(-k,h)
\]
and
\begin{align*}
\varphi(a_1) + \cdots + & \varphi(a_h) - \varphi(a'_1) - \cdots - \varphi(a'_k)\\
& = (q_1m_1+s_1,\ldots,q_dm_d+s_d) \\
& = (s_1,\ldots,s_d) + (q_1m_1,\ldots,q_dm_d) \\
& \in \varphi(hA-kA) + \Lambda_G(-k,h).
\end{align*}

(iii)  If $(s_1,\ldots,s_d) \in \varphi(hA -kA)$, then there exist elements $a_1,\ldots,a_h, a'_1,\ldots, a'_k \in G$ 
such that
\[
(s_1,\ldots,s_d) = \varphi(a_1+\cdots + a_h - a'_1 - \cdots - a'_k).
\]
Since the $i$th component of $
\varphi(a_1) + \cdots + \varphi(a_h) - \varphi(a'_1) - \cdots - \varphi(a'_k)$ is~\eqref{MSTD:lat6}, we have
\begin{align*}
\varphi(a_1) + \cdots + & \varphi(a_h) - \varphi(a'_1) - \cdots - \varphi(a'_k)\\
& = (q_1m_1+s_1,\ldots,q_dm_d+s_d) \\
& = (s_1,\ldots,s_d) + (q_1m_1,\ldots,q_dm_d).
\end{align*}
Equivalently,
\begin{align*}
(s_1,\ldots,s_d) & =
\varphi(a_1) + \cdots +  \varphi(a_h) - \varphi(a'_1) - \cdots - \varphi(a'_k) - (q_1m_1,\ldots,q_dm_d)\\
&  \in h\varphi(A) - k\varphi(A) - \Lambda_G(-k,h) \\
& = h\varphi(A) - k\varphi(A) + \Lambda_G(-h+1,k+1).
\end{align*}
This completes the proof.
\end{proof}

\bt   \label{MSTD:theorem:latineq}
Let $A$ be a nonempty subset of the finite abelian group $G = \Z/m_1\Z \times \cdots \times \Z/m_d\Z$.  Let $\varphi: G \rightarrow P_G$ be the canonical embedding defined 
by~\eqref{MSTD:map1},~\eqref{MSTD:map2}, and~\eqref{MSTD:map3}.  Let
\[
B_t = \varphi(A)+\Lambda_G(0,t) \subseteq \Zd.
\]
Then
\[
hB_t-kB_t \subseteq \varphi(hA-kA) + \Lambda_G(-kt,ht)
\]
and
\[
hB_t-kB_t \supseteq \varphi(hA-kA) 
+ \Lambda_G(-kt+h+k-1,ht-h-k+1).
\]
Moreover,
\[
|hB_t-kB_t| \leq |hA-kA| ((h+k) t)^d
\]
and
\[
|hB_t  -  kB_t| \geq  |hA-kA|((h+k)t-2(h+k-1))^d.
\]
\et

\begin{proof}
By Lemma~\ref{MSTD:lemma:lat} (ii),
\[
h\varphi(A) - k\varphi(A) \subseteq \varphi(hA-kA)  +  \Lambda_G(-k,h).
\]
It follows from~\eqref{MSTD:sumLambda}--\eqref{MSTD:sumdiffLambda} that
\begin{align*}
hB_t-kB_t & = h\left(\varphi(A)+\Lambda_G(0,t)\right) - k\left(\varphi(A)+\Lambda_G(0,t)\right) \\
& =  h\varphi(A) - k\varphi(A)+ h\Lambda_G(0,t) - k\Lambda_G(0,t) \\
& \subseteq \varphi(hA-kA)  +  \Lambda_G(-k,h) + \Lambda_G(-kt+k,ht-h+1) \\
& = \varphi(hA-kA)  +  \Lambda_G(-kt,ht).
\end{align*}
Applying~\eqref{MSTD:directsum} and~\eqref{MSTD:sizeLambda}, we obtain
\begin{align*}
|hB_t-kB_t| & \leq  |\varphi(hA-kA)  +  \Lambda_G(-kt,ht)| \\
& = |\varphi(hA-kA)| | \Lambda_G(-kt,ht)| \\
& = |hA-kA| ((h+k)t)^d.
\end{align*}

Similarly, by Lemma~\ref{MSTD:lemma:lat} (iii),
\[
\varphi(hA-kA) \subseteq 
h\varphi(A) -  k\varphi(A) + \Lambda_G(-h+1,k+1). 
\]
and so
\begin{align*}
hB_t  - kB_t & =  h\left( \varphi(A) + \Lambda_G(0,t) \right)
- k\left( \varphi(A) + \Lambda_G(0,t)\right)  \\
&  =  h\varphi(A) - k\varphi(A) + \Lambda_G(-kt+k,ht-h+1) \\
&  =  h\varphi(A) - k\varphi(A) + \Lambda_G(-h+1,k+1)+ \Lambda_G(-kt+h+k-1,ht-h-k+1) \\
& \supseteq  \varphi(hA-kA)   + \Lambda_G(-kt+h+k-1,ht-h-k+1).
\end{align*}
Therefore, 
\[
|hB_t  -  kB_t|  \geq |hA-kA|((h+k)t-2(h+k-1))^d.
\]
\end{proof}

\bt           \label{MSTD:theorem:ratio}
 Let   $G = \Z/m_1\Z \times \cdots \times \Z/m_d\Z$ be a finite abelian group, and let $\varphi$ be the canonical embedding of $G$ into the fundamental integer parallelepiped $P_G$ of the lattice $\Lambda_G$.  Let $A$ be a nonempty subset $G$ and let $h_1,h_2,k_1,k_2$ be nonnegative integers such that $h_1 \geq 1, h_2 \geq 1$, 
\[
h_1+k_1 = h_2+k_2
\]
and
\[
|h_1A-k_1A| > |h_2A - k_2A|.
\]
For every sufficiently large integer $t$, the set of lattice points
\[
B_t = \varphi(A) + \Lambda_G(0,t)
\]
satisfies the inequality
\[
|h_1B_t-k_1B_t| > |h_2B_t - k_2B_t|.
\]
\et
 
\begin{proof}
Let $h_1+k_1 = h_2+k_2 = c.$
By Theorem~\ref{MSTD:theorem:latineq}, 
\begin{align*}
\frac{|h_1B_t-k_1B_t|}{|h_2B_t - k_2B_t|}
& \geq \frac{|h_1A-k_1A|((h_1+k_1)t-2(h_2+k_2 -1))^d}{|h_2A-k_2A| ((h_2+k_2) t)^d} \\
& \geq \frac{|h_1A-k_1A|}{|h_2A-k_2A|}\frac{(ct-2c+2)^d}{(ct)^d} \\
& = \frac{|h_1A-k_1A|}{|h_2A-k_2A|}
\left( 1 - \frac{2c-2}{c t}\right)^d
\end{align*}
Since
\[
\lim_{t\rightarrow \infty} \left( 1 - \frac{2c-2)}{c t}\right)^d =1
\]
it follows that 
\[
|h_1B_t-k_1B_t| > |h_2B_t - k_2B_t|
\]
for all sufficiently large $t$.  This completes the proof.
\end{proof}

\bt           \label{MSTD:theorem:sumdiffratio}
 Let   $G = \Z/m_1\Z \times \cdots \times \Z/m_d\Z$ be a finite abelian group, and let $\varphi$ be the canonical embedding of $G$ into the fundamental integer parallelepiped $P_G$ of the lattice $\Lambda_G$.  Let $A$ be a nonempty subset $G$ such that
 \[
|A+A| > |A - A|.
\]
For every sufficiently large integer $t$, the set of lattice points
\[
B_t = \varphi(A) + \Lambda_G(0,t)
\]
satisfies the inequality
\[
|B_t + B_t| > |B_t - B_t|.
\]
\et

\begin{proof}
Apply Theorem~\ref{MSTD:theorem:ratio} with $(h_1,k_1) = (2,0)$ and 
$(h_2,k_2) = (1,1).$
\end{proof} 

The final step is to embed a finite set of lattice points in \Zd\ into the integers.  For $a = (a_1,\ldots, a_d) \in \Zd$, define
\[
\| a\| = \max( |a_i| : i=1,\ldots,d).
\]

\bl  \label{MSTD:lemma:psi}
Let $m$ be a positive integer, and consider the group homomorphism
$\psi:\Zd \rightarrow \Z$ defined by 
\[
\psi(a) = \sum_{i=1}^d a_i m^{i-1}
\]
for all $a = (a_1,\ldots,a_d) \in \Zd.$
If $\| a\| < m$ and $\psi(a) = 0$, then $a=0.$
\el

\begin{proof}
Suppose that $a \in \Zd$ satisfies $\| a\| < m$ and $\psi(a) = 0$.
The inequality  $\| a\| < m$ implies that $|a_i| 
\leq m-1$ for all $i=1,\ldots,d$.
If $a \neq 0$, let $r$ be the largest integer such that $a_r \neq 0$.
Then $\psi(a)=0$ implies that
\[
a_rm^{r-1} = -\sum_{i=1}^{r-1} a_im^{i-1}
\]
but this is impossible since 
\[
|a_rm^{r-1} | \geq m^{r-1}
\]
and
\[
\left| -\sum_{i=1}^{r-1} a_im^{i-1} \right|
\leq m^{r-1}-1.
\]
This completes the proof.
\end{proof}

\bt
Let $L$ be a positive integer, and let $A$ be a nonempty finite subset of \Zd.  
Let $m$ be an integer such that
\[
m > 2L\max( \| a\| : a \in A).
\]
Define $\psi: \Zd \rightarrow \Z$ by
\[
\psi(a_1,\ldots,a_d) = \sum_{i=1}^d a_i m^{i-1}.
\]
If $h$ and $k$ are nonnegative integers such that 
\[
h+k \leq L
\]
then
\[
|hA-kA| = |h\psi(A) - k\psi(A)|.
\]
\et

\begin{proof}
Since $\psi$ is a homomorphism, we have
\[
\psi(hA-kA) = h\psi(A)-k\psi(A).
\]
We must show that the map $\psi: hA-kA \rightarrow \Z$ is one-to-one.  Let $u,v \in hA-kA$ such that $\psi(u) = \psi(v)$.  The set $A$ contains lattice points
$a_j = (a_{1,j},\ldots,a_{d,j})$ and $b_j = (b_{1,j},\ldots,b_{d,j})$ for $j = 1,\ldots,h$, and lattice points
$a'_{\ell} = (a'_{1,{\ell}},\ldots,a'_{{d,\ell}})$ and $b'_{\ell} = (b'_{1,{\ell}},\ldots,b'_{d,{\ell}})$ for ${\ell} = 1,\ldots,k$  such that 
\[
u = a_1 + \cdots + a_h - a'_1 - \cdots - a'_k
\]
and
\[
v = b_1 + \cdots + b_h - b'_1 - \cdots - b'_k.
\]
The $i$th coordinate of $u$ is $\sum_{j=1}^h a_{i,j} - \sum_{\ell = 1}^k a'_{i,\ell}$, and so
\begin{align*}
\| u\| & \leq \sum_{j=1}^h \|a_j\| + \sum_{\ell =1}^k \|a'_k\| \\
& \leq (h+k)\max( \| a\| : a \in A) \\
& < \frac{m}{2}.
\end{align*}
Similarly,
\[
\| v\| \leq \sum_{j=1}^h \|b_j\| + \sum_{\ell =1}^k \|b'_k\| < \frac{m}{2}
\]
and so
\[
\| u-v\| < m.
\]
If $\psi(u) = \psi(v)$, then $\psi(u-v) = 0$.  
By Lemma~\ref{MSTD:lemma:psi}, it follows that $u=v$ and the map $\psi$ is one-to-one.
\end{proof}

\section{A counting argument}
In this section we prove that the finite abelian group
\[
G = \Z/n\Z \times \Z/2\Z
\]
contains MSTD sets.

Let $\Omega$ be the set of all subsets of $G$ of the form
\[
A = \{ (i+n\Z,\varepsilon_i+2\Z) : i = 0,1,\ldots,n-1 \text{ and } \varepsilon_i\in \{0,1\} \}.
\]
Then $|\Omega| = 2^n$ and $|A|=n$ for all $A \in \Omega$.
If $g  =  (b+n\Z,\delta+2\Z)$ is in the difference set $A-A$, then there exist
$a = (i+n\Z,\varepsilon_i+2\Z)\in A$ and $a' = (j+n\Z,\varepsilon_j+2\Z)\in A$ such that $g = a-a'$, or, equivalently,
\begin{align*}
(b+n\Z,\delta+2\Z) & = (i+n\Z,\varepsilon_i+2\Z) - (j+n\Z,\varepsilon_j+2\Z)\\ &= (i-j+n\Z,\varepsilon_i-\varepsilon_j+2\Z).
\end{align*}
If $b+n\Z =  n\Z$, then $i\equiv j\pmod{n}$ and so $i=j$.  It follows  that $\varepsilon_i = \varepsilon_j$ and  $\delta\equiv 0\pmod{2}$.
This implies that for every set $A \in \Omega$,  if $(n\Z,\delta+2\Z) \in A-A$, then $\delta\equiv 0 \pmod{2}$, that is,
\[
(n\Z,1+2\Z) \notin A-A
\]
and so
\[
|A-A| \leq |G|-1 = 2n-1.
\]
We shall prove that there exists $A \in \Omega$ such that $A+A=G$, that is,
\[
|A+A|= 2n
\]
and so $A$ is an MSTD set.
Indeed, we shall prove that $A+A=G$ for almost all sets $A\in \Omega$.

\bt
Let
\[
G = \Z/n\Z \times \Z/2\Z
\]
and let $\Omega$ be the set of all subsets of $G$ of the form
\[
A = \{ (i+n\Z,\varepsilon_i+2\Z) : i = 0,1,\ldots,n-1 \text{ and } \varepsilon_i\in \{0,1\} \}.
\]
Let $\Psi(G)$ denote the number of $A\in \Omega$ such that $A+A=G.$  Then
\[
\Psi(G) \geq 
\begin{cases}
2^n \left( 1 - \frac{\sqrt{2}n}{2^{n/2}}  \right) & \text{if $n$ is odd} \\ 
2^n \left( 1 - \frac{2n}{2^{n/2}}  \right) & \text{if $n$ is even.} 
\end{cases}
\]
\et 

\begin{proof}
Let 
\[
g  =  (b+n\Z,\delta+2\Z)\in G
\]
and let $\phi(g)$ denote the number of $A \in \Omega$ such that $g\notin A+A$.  Then
\[
|\Omega| - \Psi(G) \leq \sum_{g\in G} \phi(g).
\]

Consider first the case that $n$ is odd.  There is a unique congruence class $a_0+n\Z \in \Z/n\Z$ such that $2a_0 \equiv b\pmod{n}$, and $(a_0 +n\Z,\varepsilon+2\Z)\in A$ for some $\varepsilon \in \{0,1\}.$
If $\delta \equiv 0 \pmod{2}$, then 
\[
g = (b+n\Z,2\Z) = (a_0 +n\Z,\varepsilon+2\Z) + (a_0 +n\Z,\varepsilon+2\Z) \in A+A
\]
and $\phi(g) = 0$.

If $\delta \equiv 1 \pmod{2}$, then $(a_0+n\Z,\varepsilon+2\Z) + (a_0+n\Z,\varepsilon+2\Z) \neq g$.
If $a\not\equiv a_0\pmod{n}$, then $b-a\not\equiv a\pmod{n}$.  
It follows that
$\Z/n\Z \setminus \{a_0+n\Z\}$ is the union of $(n-1)/2$ disjoint sets of the form $\{a_j+n\Z,b-a_j+n\Z\}$, for $j=1,\ldots,(n-1)/2$.
Let $\{\varepsilon_j\}_{j=0}^{(n-1)/2}$ be any sequence of 0's and 1's.  We define the sequence $\{\varepsilon'_j\}_{j=1}^{(n-1)/2}$ by
\[
\varepsilon'_j = \delta-\varepsilon_j+1.
\]
Defining the set $A \in \Omega$ by
\[
A = \{ (a_j+n\Z,\varepsilon_j+2\Z\}_{j=0}^{(n-1)/2}
\cup \{ (b-a_j+n\Z,\varepsilon'_j+2\Z\}_{j=1}^{(n-1)/2}
\]
we see that $g \notin A+A$.
It follows that
\[
\phi(g) = 
\begin{cases}
0 & \text{if $\delta \equiv 0 \pmod{2}$} \\
 2^{(n+1)/2} & \text{if $\delta \equiv 1 \pmod{2}$} \\
\end{cases}
\]
and so
\[
\sum_{g\in G} \phi(g) = \left( \frac{|G|}{2}\right) 2^{(n+1)/2} = n2^{(n+1)/2}
\]
and
\begin{align*}
\Psi(G) & \geq |\Omega| - \sum_{g\in G} \phi(g) \\
& = 2^n - n2^{(n+1)/2} \\
& = 2^n\left( 1 - \frac{\sqrt{2} n}{2^{n/2}} \right).
\end{align*}

There is a similar argument in the case that $n$ is even. 
If $b$ is odd, then the congruence $2x\equiv b \pmod{n}$ has no solution, and $\Z/n\Z$ can be partitioned into $n/2$ pairs of elements that sum to $b$.  It follows that $\phi(g) = 2^{n/2}$.  There are exactly $n/2$ elements in the group $\Z/n\Z$ with odd $b$, and so
$n$ elements in the group $G$ with odd $b$.
Therefore,
\[
\sum_{  \substack{g\in G \\  b \text{ odd } }  } \phi(g) = n2^{n/2}.
\]

If $b$ is even, then there is an integer $a_0$ such that 
\[
b+n\Z = 2(a_0+n\Z) = 2(a_0+n/2 + n\Z).
\]
If $\delta \equiv 0\pmod{2}$, then $\phi(g)=0$.
If $\delta  \equiv 1\pmod{2}$, then $\Z/n\Z \setminus \{a_0+n\Z,a_0+n/2 + n\Z\}$ can be partitioned into $(n-2)/2$ pairs of elements that sums to $b$.  It follows that $\phi(g) = 4\cdot 2^{(n-2)/2} = 2^{(n+2)/2}$ and 
\[
\sum_{  \substack{g\in G \\  b \text{ even} }  } \phi(g) = \left(\frac{n}{2}\right) 2^{(n+2)/2} = n2^{n/2}.
\]
Therefore,
\[
\sum_{  g\in G  } \phi(g) =  2n2^{n/2}
\]
and
\begin{align*}
\Psi(G) & \geq |\Omega| - \sum_{g\in G} \phi(g) \\
& = 2^n - n2^{(n+2)/2} \\
& = 2^n\left( 1 - \frac{2n}{2^{n/2}} \right).
\end{align*}
This completes the proof.
\end{proof}

\section{Problems}
It would be interesting to classify the structure of all MSTD sets of integers.  We can stratify this problem in the following way.  
Let $\mathcal{F}(\Z)$ denote the set of all finite sets of integers.  Define the function $n: \mathcal{F}(\Z) \rightarrow \Z$ by
\[
n(A) = |A+A|-|A-A|.
\]
What is the range of this function?  Is it possible to construct, for every integer $t$, a set $A \in \mathcal{F}(\Z)$ with $n(A) = t$?
Can we describe all sets $A$ with $n(A) = 1$?

{\it Acknowledgements.}  The idea of constructing of MSTD sets of lattice points and integers from MSTD subsets of finite abelian groups is due to Terence Tao~\cite{tao06x}, and I thank him for his permission to include it in this paper.  I also thank Boris Bukh and Peter Hegarty for helpful discussions.

\providecommand{\bysame}{\leavevmode\hbox to3em{\hrulefill}\thinspace}
\providecommand{\MR}{\relax\ifhmode\unskip\space\fi MR }
\providecommand{\MRhref}[2]{%
  \href{http://www.ams.org/mathscinet-getitem?mr=#1}{#2}
}
\providecommand{\href}[2]{#2}

\end{document}